\documentclass{article}%
\usepackage{amsfonts}
\usepackage{amsmath}
\usepackage{amssymb}
\usepackage{graphicx}%
\newtheorem{theorem}{Theorem}

\newtheorem{corollary}{Corollary}

\newtheorem{definition}{Definition}
\newtheorem{example}{Example}
\newtheorem{lemma}{Lemma}
\newtheorem{notation}{Notation}
\newtheorem{proposition}{Proposition}
\newtheorem{remark}{Remark}
\newenvironment{proof}[1][Proof]{\noindent\textbf{#1.} }{\ \rule{0.5em}{0.5em}}
\begin{document}

\title{The Automorphism Group of Certain Factorial Threefolds and a Cancellation Problem}
\author{David R. Finston and Stefan Maubach\\David R. Finston\\Department of Mathematical Sciences\\New Mexico State University\\Las Cruces, New Mexico 88003\\dfinston@nmsu.edu\\Stefan Maubach\\Department of Mathematical Sciences\\University of Texas at Brownsville\\Brownsville, TX 78520\\stefan.maubach@utb.edu}
\date{}
\maketitle

\begin{abstract}
The automorphism groups of certain factorial complex affine threefolds
admitting locally trivial actions of the additive group are determined. \ As a
consequence new counterexamples to a generalized cancellation problem are obtained.

\end{abstract}

\section{Introduction}

A well known cancellation problem asks, for complex affine varieties $X$ and
$Y,$ whether an isomorphism $X\times\mathbb{C\cong}Y\times\mathbb{C}$ implies
that $X\cong Y.$ For $X$ and $Y$ of dimension 1 a positive answer is given by
\cite{AEH} and for $X$ and $Y$ of dimension 2 counterexamples are provided by
the Danielewski surfaces \cite{D} \cite{ML} \cite{Fi} \cite{tD}. \ On the
other hand, for $X\times\mathbb{C}\cong\mathbb{C}^{3},$ Fujita and
Miyanishi-Sugie proved that $X\cong\mathbb{C}^{2}.$ \ The Danielewski surfaces
can be realized as total spaces for principal bundles for $G_{a},$ the
additive group of complex numbers, over the affine line with two origins. They
are therefore smooth surfaces, but nonfactorial, i.e. their coordinate rings
lack the unique factorization property. It is natural then to ask whether the
cancellation problem has a positive solution for factorial affine varieties,
or for affine total spaces of principal $G_{a}$ bundles over quasiaffine
varieties. \ We produce families of three dimensional counterexamples.

To point out the role played by principal $G_{a}$ bundles, let $Y\ $be a
scheme over $\mathbb{C},$ and $X_{1},X_{2}$ total spaces for principal $G_{a}$
bundles over $Y.$ Then each $X_{i}$ is represented by a one cocycle in
$H^{1}(Y,O_{Y}),$ and we can represent the base extension $X_{1}\times
_{Y}X_{2}$ by elements of $H^{1}(X_{1},O_{X_{1}})$ (and $H^{1}(X_{2},O_{X_{2}%
}).$ \ If the $X_{i}$ are affine then $H^{1}(X_{i},O_{X_{i}})=0$ and therefore
$X_{1}\times\mathbb{C}\cong X_{1}\times_{Y}X_{2}\cong X_{2}\times\mathbb{C}$.
\ In particular, affine total spaces for principal $G_{a}$ bundles is a
natural context in which to seek potential counterexamples to the cancellation problem.

In the case of the Danielewski surfaces, not only are the bundles
inequivalent, the total spaces are not homeomorphic in the natural (complex)
topology on $\mathbb{C}^{3}$, let alone isomorphic as varieties. \ For a
complex quasiprojective base however, a principal $G_{a}$ bundle is
necessarily trivial in the natural topology \cite{Steenrod}. \ Thus algebraic
methods are necessary to distinguish the total spaces. The Makar-Limanov
invariant, which for an affine $k-$domain $A$ is the intersection of the
kernels of all locally nilpotent $k-$derivations of $A$, provides the
necessary algebraic tool enabling the determination of the automorphism groups
of certain affine threefolds, all obtained as total spaces for principal
$G_{a}$ bundles over the spectrum of singular but factorial complex surfaces
punctured at the singular point. \ A class of these threefolds yield the
desired counterexamples:

\begin{example}
Let $X_{n,m}\subset\mathbb{C}^{5}$ be the affine variety defined by
\[
X^{a}+Y^{b}+Z^{c},X^{m}U-Y^{n}V-1
\]
with $m,n$ positive integers and $a,b,c\ $pairwise relatively prime positive
integers satisfying $\frac{1}{a}+\frac{1}{b}+\frac{1}{c}<1.$ Then $X_{n,m}$ is
factorial,
\[
X_{n,m}\times\mathbb{C}\cong X_{n^{\prime},m^{\prime}}\times\mathbb{C}%
\]
$\mathbb{\ }$for all $(m,n),(m^{\prime},n^{\prime}),$ but $X_{n,m}\ \cong
X_{n^{\prime},m^{\prime}}\ $implies that $(m,n)=(m^{\prime},n^{\prime}).$
\end{example}

We suspect that the condition $\frac{1}{a}+\frac{1}{b}+\frac{1}{c}<1$ can be weakened.

Principal $G_{a}$ bundles with affine total space $X$ arise from locally
trivial algebraic $G_{a}$ actions on $X$. \ The local triviality implies that
the quotient $X/G_{a}$ exists as an algebraic scheme, and gives $X$ the
structure of a principal $G_{a}$ bundle over $X/G_{a}.\ $ \ If $X$ is in
addition factorial, then $X/G_{a}$ has the structure of a quasiaffine variety.
\ The Makar-Limanov invariant enters the picture since every algebraic $G_{a}$
action on an affine $X$ arises as the exponential of a locally nilpotent
derivation $D$ of $\mathbb{C}[X].$ If $X$ is factorial, then the action is
locally trivial if and only if $\ker(D)\cap im(D)$ generates the unit ideal in
$\mathbb{C}[X]$ \cite{D-F2}. \

\section{The Makar-Limanov Invariant.}

The condition on the exponents $a,b,c$ in the above example will enable us to
use Mason's theorem, stated here as Theorem 1. { Let $k$ be a field of
characteristic }${0}$ { and, for $f\in k[T]$, denote by $N(f)$ the number of
distinct zeroes of $f$ in an algebraic closure of }$k${ . }

\begin{theorem}
(e.g. \cite{Snyder}) \label{Mason} Let $f,g\in k[T]$ and let $h=f+g$. Assume
that $f,g,h$ are relatively prime of positive degree. Then
\[
{max\{deg(f),deg(g),deg(h)\}<N(fgh).\ }%
\]

\end{theorem}

Two corollaries apply to the problem at hand.

\begin{corollary}
{ \label{05.10.cor1} Let $P(X,Y,Z)=X^{a}+Y^{b}+Z^{c}$ where $a,b,c\in
\mathbb{N}$ satisfy}%
\[
{ \ \frac{1}{a}+\frac{1}{b}+\frac{1}{c}\leq1.}%
\]
{ If $f,g,h\in k[T]$ satisfy }
\end{corollary}

\begin{enumerate}
\item { $P(f,g,h)=0$ and }

\item { $f,g,h$ are relatively prime.}
\end{enumerate}

Then { at least one of $f,g,h$ must be constant. }

\begin{proof}
{ It is enough to consider the case that $k$ is algebraically closed.
\ }Assume that none of $f,g,h$ is constant. {Applying Mason's theorem, \ the
fact that $f^{a}+g^{b}+h^{c}=0$ yields:}%
\[%
\begin{tabular}
[c]{lll}%
${ max(a\cdot deg(f),b\cdot deg(g),c\cdot deg(h))}$ & $<$ & ${ N(f^{a}%
)+N(g^{b})+N(h^{c})}$\\
& $=$ & ${ N(f)+N(g)+N(h)}$\\
& $\text{\ }{ \leq}$ & ${ deg(f)+deg(g)+deg(h).}$%
\end{tabular}
\]

{ Suppose }%
\[
{ a\cdot deg(f)\geq b\cdot deg(g),\ a\cdot deg(f)\geq c\cdot deg(h)>0.}%
\]
{ Then }$\frac{\deg g}{\deg f}\leq\frac{a}{b},\frac{\deg h}{\deg f}\leq
\frac{a}{c}$ so that%
\[%
\begin{tabular}
[c]{lll}%
${ a\cdot deg(f)}$ & ${ =}$ & ${ max(a\cdot deg(f),b\cdot deg(g),c\cdot
deg(h))}$\\
& $<$ & $deg(f)+deg(g)+deg(h)$\\
& $\leq$ & $deg(f)(1+\frac{a}{b}+\frac{a}{c}).$%
\end{tabular}
\]
{ Thus}${ \ 1<\frac{1}{a}+\frac{1}{b}+\frac{1}{c}}${ , which exactly
contradicts the assumption. }

{ The cases where $b\cdot deg(g)$ or $c\cdot deg(h)$ is the largest go
equivalently. }
\end{proof}

\begin{corollary}
{ \label{05.10.cor2} Let $P(X,Y,Z)=X^{a}+Y^{b}+Z^{c}+\lambda$ where
$\lambda\in k$, and $a,b,c\in\mathbb{N}\backslash\{0,1,2,3\}$ satisfy
$\frac{1}{a-3}+\frac{1}{b-3}+\frac{1}{c-3}\leq\frac{1}{2}$. If $f,g,h\in k[T]$
satisfy }
\end{corollary}

\begin{enumerate}
\item { $P(f,g,h)=0$ and}

\item $f,g,h$ are relatively prime.
\end{enumerate}

Then {at least one of $f,g,h$ must be constant. }\newline

\begin{proof}
Again\emph{ }{ it is enough to consider the case that $k$ is algebraically
closed. We will arrive at a contradiction from the assumption that }%
$f^{a}+g^{b}+h^{c}=\lambda$ for some nonconstant $f,g,h$. Taking derivatives
with respect to $T$ yields $af^{a-1}f^{\prime}+bg^{b-1}g^{\prime}%
+ch^{c-1}h^{\prime}=0$. Now we cannot apply Mason's theorem directly as there
may be common factors in $ff^{\prime},gg^{\prime},hh^{\prime}$. Define
$w:=gcd(f^{a-1}f^{\prime},g^{b-1}g^{\prime},h^{c-1}h^{\prime})$. Using the
fact that $gcd(xy,z)$ divides $gcd(x,z)gcd(y,z)$ repeatedly we see that $w$
divides $gcd(f^{\prime},g^{b-1}g^{\prime},h^{c-1}h^{\prime})\cdot
gcd(f^{a-1},g^{\prime},h^{c-1}h^{\prime})\cdot gcd(f^{a-1},g^{b-1},h^{\prime
})\cdot gcd(f^{a-1},g^{b-1},h^{c-1})$ and since $gcd(f,g,h)=1$, we see that
$deg(w)\leq deg(f^{\prime})+deg(g^{\prime})+deg(h^{\prime}%
)=deg(f)+deg(g)+deg(h)-3$. One can apply Mason's theorem to
\[
a\frac{1}{w}f^{a-1}f^{\prime}+b\frac{1}{w}g^{b-1}g^{\prime}+c\frac{1}%
{w}h^{c-1}h^{\prime}=0,
\]
which, together with some calculus, yields
\[%
\begin{array}
[c]{rl}
& 2(deg(f)+deg(g)+deg(h))\\
\geq & N(ff^{\prime}gg^{\prime}hh^{\prime})\\
\geq & N(ff^{\prime}\frac{1}{w}gg^{\prime}\frac{1}{w}hh^{\prime}\frac{1}%
{w})\\
\text{(Mason's)}>&max\big(deg(f^{a-1}f^{\prime}\frac{1}{w}),deg(g^{b-1}g^{\prime}\frac{1}%
{w}),deg(h^{c-1}h^{\prime}\frac{1}{w})\big)\\
= & max\big(deg(f^{a-1}f^{\prime}),deg(g^{b-1}g^{\prime}\big),deg(h^{c-1}%
h^{\prime})\big)-deg(w)\\
\geq&max\big(deg(f^{a-1}f^{\prime}),deg(g^{b-1}g^{\prime}),deg(h^{c-1}%
h^{\prime})\big)\\
&~~~~~~~-deg(f)-deg(g)-deg(h)+3\\
=& max\big(a deg(f)-1,b deg(g)-1,c deg(h)-1\big)\\
& ~~~~~~~-deg(f)-deg(g)-deg(h)+3\\
\geq & max\big((a-3)deg(f),(b-3)deg(g),(c-3)deg(h)\big)+2\\
> & max\big((a-3)deg(f),(b-3)deg(g),(c-3)deg(h)\big)
\end{array}
\]
Assuming that $max((a-3)deg(f),(b-3)deg(g),(c-3)deg(h))=(a-3)deg(f)$ (the
other cases go similarly) then will yield $(a-3)deg(f)<2(1+\frac{a-3}%
{b-3}+\frac{a-3}{c-3})deg(f)$ which exactly contradicts the assumption
$\frac{1}{a-3}+\frac{1}{b-3}+\frac{1}{c-3}\leq\frac{1}{2}$.
\end{proof}

\begin{definition}
\begin{enumerate}
\item For a $k$-domain $B$, $LND(B)$ is the set of locally nilpotent $k$
derivations of $B.$ \

\item Given $D\in LND(B)$, $s\in B$ is a slice for $D$ if $D(s)=1.$

\item Given $D\in LND(B),$ an element $p$ of $B$ is called a preslice if
$0=D^{2}(p)\neq D(p).$
\end{enumerate}
\end{definition}

\begin{remark}
A preslice always exists for a nonzero locally nilpotent derivation $D.$
\ Indeed, by local nilpotency, for $b\in B\ -\ker(D)$, there is a positive
integer $n$ for which $0\ \neq D^{n+1}(b)\in\ker(D).$ \ Then $p=D^{n}(b)$ is a
preslice. \ If $D$ admits a slice $s$, then $B=B^{D}[s]$, where $B^{D}$
denotes $\ker(D),$ and therefore $D$ $=$ $\frac{\partial}{\partial s}$
\cite{D-F}$.$
\end{remark}

\begin{lemma}
{ \label{05.10.lemma1} Let $A$ be a $\mathbb{C}$-domain and $x,y,z\in
A\backslash\{0\}$. Let $P=x^{a}+y^{b}+z^{c}+\lambda$ for some $a,b,c\in
\mathbb{N}\backslash\{0,1\}$, $\lambda\in\mathbb{C}$. Let $B:=A/(P)$, and
assume that $B$ is a domain (i.e. $P$ is a prime element of }$A${ ). If either
\newline i) $\lambda=0$ and $\frac{1}{a}+\frac{1}{b}+\frac{1}{c}\leq1$, or
\newline ii) $a,b,c\geq4$ and $\frac{1}{a-3}+\frac{1}{b-3}+\frac{1}{c-3}%
\leq\frac{1}{2}$, \newline then $D\in LND(B)$ implies $D(x)=D(y)=D(z)=0$. }
\end{lemma}

\begin{proof}
Since { $B$ is a domain, and $D$ is locally nilpotent, a preslice $p$ exists.
Set $q:=D(p)$ (and thus $q\in B^{D}$) and observe that }$D$ extends uniquely
to a locally nilpotent derivation { $\tilde{D}$ of $\tilde{B}:=B[q^{-1}]$.
\ Since $\tilde{D}$ has the slice $s:=p/q$ we have $\tilde{B}=\tilde{B}%
^{D}[s]$. We can identify $\tilde{D}$ with $\frac{\partial}{\partial s}$.
Denote by $K\ $\ the quotient field of $\tilde{B}^{\frac{\partial}{\partial
s}}$ ($=$ quotient field of $B^{D})$ noting that }$D$ { extends uniquely to }
$\frac{\partial}{\partial s}${ on $K[s]$. \ Write $x,y,z\in K[s]$, as
$x=f(s),y=g(s),z=h(s)$ for some polynomials $f,g,h\in K[s]$. If }%
${k=\gcd(f,g)}$ then $k$ divides $h$ as well. \ Writing
\[
f=k\widehat{f},g=k\widehat{g},h=k\widehat{h}%
\]
we obtain%
\[
(k^{bc}\widehat{f})^{a}+(k^{ac}\widehat{g})^{b}+(k^{ab}\widehat{h})^{c}=0
\]
and therefore%
\[
\widehat{f}^{a}+\widehat{g}^{b}+\widehat{h}^{c}=0
\]
with $\widehat{f},\widehat{g},\widehat{h}$ pairwise relatively prime. \

{ In case i) we can use corollary \ref{05.10.cor1}, to conclude that }$k$ and
at least one of $\widehat{f},\widehat{g},\widehat{h}$ lie in $K$, so that one
of $x,y,z$ lies in $\ker(D).$ But if, for instance, $D(x)=0,$ then
$0=D(y^{b}+z^{c})$ then { by the following lemma we see that $D(y)=D(z)=0$. }

{ Similarly in case ii) we can use corollary \ref{05.10.cor2} to conclude that
at least one of $x,y,z$ must lie in }$\ker(D)${ . Suppose it is $x$. Then
again $D(y^{b}+z^{c})=0$ where $b,c\geq2$. }
\end{proof}

\begin{lemma}
{ \label{lenny} (Makar-Limanov \cite[Lemma 2]{lenny} ) Let $A$ be a domain and
let $n,m\in\mathbb{N}$ satisfying $n,m\geq2$. If $D\in LND(A)$ and
$D(c_{1}a^{n}+c_{2}b^{m})=0$ where $a,b\in A$, $c_{1},c_{2}\in A^{D}$, and
$c_{1}a^{n}+c_{2}b^{m}\neq0$ . Then $D(a)=D(b)=0$. }
\end{lemma}

{ Fix $P(X,Y,Z):=X^{a}+Y^{b}+Z^{c}+\lambda$ in }$\mathbb{C}[X,Y,Z]$ { and
assume that $P$ is irreducible, i.e. that }$a,b,c$ are pairwise relatively
prime.{ }

\begin{notation}
{For the remainder of the paper, $R:=\mathbb{C}[X,Y,Z]/(P)$, and $x,y,z$
denote the images of $X,Y,Z$ in $R.$ Set $A_{n,m}:=R[U,V]/(x^{m}U-y^{n}V-1)$
where $m,n\in\mathbb{N}$, $m,n\geq2$. The images of $U,V$ in $A_{n,m}$ will be
denoted by }$u,v.$
\end{notation}

\begin{proposition}
{ \label{05.10.UFD} If $gcd(a,b)=gcd(a,c)=gcd(b,c)=1$, then $A_{n,m}$ is a
UFD. }
\end{proposition}

\begin{proof}
That $R$ is a UFD in case { $\lambda=0$ is a well known result of Samuel. \ A
slight modification of the argument in \cite{Samuel} yields the result for
$\lambda\neq0.$ Define an }$R$ derivation $D$ of $A_{n,m\text{ }}$ by setting
$D(v)=x^{m},$ $D(u)=y^{n}).$ Clearly $D$ is { locally nilpotent and generates
a locally trivial }$G_{a}$ action on the smooth variety $X_{n,m}\equiv$
\textbf{Spec }$A_{n,m\text{ }}$. \ The quotient $X_{n,m}/G_{a}$ is isomorphic
to the complement of a finite but nonempty subset of \textbf{Spec} $R.$ \ The
quotient map $X_{n,m}\rightarrow X_{n,m}/G_{a}$ is a Zariski fibration with
both the base and fiber having trivial Picard group. \ By \cite{Magid} we
conclude that $Pic(X_{n,m})$ is also trivial and therefore $A_{n,m\text{ }}$is
a UFD.

In case { $\lambda=0$ one can argue directly that }$A_{n,m\text{ }}$is a UFD
{using Nagata's theorem \cite[Theorem 20.2]{Mat} . \ Note that }$x$ is a prime
element in $A_{n,m}:$%
\begin{align*}
A_{n,m}/(x) &  \cong\mathbb{C}[Y,Z,U,V]/(Y^{n}V+1)\\
&  \cong\mathbb{C}[Y,Z,U][\frac{1}{Y}]
\end{align*}
a domain, and
\[
A_{n,m\text{ }}[x^{-1}]\cong\mathbb{C}[X,Y,Z]/(X^{a}+Y^{b}+Z^{c})[x^{-1}][U]
\]
is a UFD.
\end{proof}

{ The following is a consequence of Lemma 1. \ref{05.10.lemma1} }

\begin{corollary}
{ \label{05.10.cor3} If $D\in LND(A_{n,m})$ then $D(x)=D(y)=D(z)=0 $. }
\end{corollary}

\begin{lemma}
{ \label{05.10.lemma2} Let $D \in LND(A_{n,m})$ and assume $D\not = 0$. Then
$A_{n,m}^{D}=\mathbb{C}[x,y,z]$. }
\end{lemma}

\begin{proof}
{ $x^{m}D(u)-y^{n}D(v)=D(x^{m}u-y^{n}v)=D(1)=0$. Since $A_{n,m}$ is a UFD we
see that $D(u)=cy^{n}$ for some $c$. Thus $D(v)=x^{m}c$. Thus $D$ is
equivalent to the locally nilpotent derivation $D^{\prime}=y^{n}\partial
_{u}+x^{m}\partial_{v}$ in particular they have the same kernel. An easy
application of the algorithm in \cite{arno} reveals that }$\ker(D^{\prime})=${
$\mathbb{C}[x,y,z]$. }
\end{proof}

\begin{theorem}
$ML${$(A_{n,m})=R.$}
\end{theorem}

\section{The Automorphism Group}

In this section we take {$R:=\mathbb{C}[X,Y,Z]/($}$X^{a}+Y^{b}+Z^{c})$ with
$a,b,c$ pairwise relatively prime satisfying
\[
{\ \frac{1}{a}+\frac{1}{b}+\frac{1}{c}<1,}%
\]
and $A_{n,m}$ as before.{ \ The derivation }%
\[
{E:=y^{n}\partial_{u}+x^{m}\partial_{v}\in Der_{\mathbb{C}}(A_{n,m}).\ }%
\]
plays a special role.

\begin{lemma}
{ \label{05.10.lemma3} Let }$B$ { be a }$k$-domain{ , and $\varphi\in Aut(B)$.
Then $\varphi^{-1}LND(B)\varphi=LND(B)$. Also, $\varphi(ML(B))=ML(B)$. }
\end{lemma}

\begin{proof}
{ If $D$ is LND, then $\varphi^{-1}D\varphi$ is also LND. So $\varphi
^{-1}LND(B)\varphi\subseteq LND(B)$ for any automorphism $\varphi$. Then }%
\[
{ \varphi^{-1}(\varphi LND(B)\varphi^{-1})\varphi\subseteq\varphi
^{-1}LND(B)\varphi,}%
\]
{ which proves the converse inclusion. }

It follows moreover that%
\[
{\varphi(ML(B))=\varphi\Big(\bigcap_{D\in LND(B)}\ker(D)\Big)}%
\]%
\[
{=\bigcap_{D\in LND(B)}\varphi(B^{D})}%
\]%
\[
=\bigcap_{D\in LND(B)}B^{\varphi D\varphi^{-1}}%
\]

{ which is equal to $ML(B)$ since $\varphi LND(B)\varphi^{-1}=LND(B)$. }
\end{proof}

\begin{corollary}
{ \label{05.10.cor4} Let $\varphi\in Aut_{\mathbb{C}}(A_{n,m})$. Then
$\varphi^{-1} E \varphi= \lambda E$ where $\lambda\in\mathbb{C}^{*}$. }
\end{corollary}

\begin{proof}
{ $LND(A_{n,m})=\mathbb{C}[x,y,z]E$, so by Lemma \ref{05.10.lemma3}
$\varphi(E)=\lambda E$ for some $\lambda\in\mathbb{C}[x,y,z]^{\ast}%
=\mathbb{C}^{\ast}$. }
\end{proof}

{ Let $S\subset T\subset B$ be domains, }$T$ { an $S$-algebra, and }$B$ { a
}$T${ -algebra. Suppose that for any $\varphi\in Aut_{S}B$ we have
$\varphi(T)=T$. Then restriction to $T$ defines a group homomorphism }$\rho
:${$Aut_{S}B\rightarrow$ $Aut_{S}T$ and $Aut_{S}B$ is an extension of
$Aut_{T}B$ by the image of }$\rho${. For }$S=\mathbb{C},T=R,B=A_{n,m}$ we will
show that $\rho$ is surjective, and determine $Aut_{\mathbb{C}}R$ and
$Aut_{R}A_{n,m}.$

\begin{remark}
Incidentally, once can see easily that no two dimensional UFD can give rise to
a counterexample to generalized cancellation via non trivial $G_{a}$ bundles.
\end{remark}

\begin{lemma}
Let A be a two dimensional finitely generated $\mathbb{C}$ algebra which is a
UFD. \ If $A$ admits a nonzero LND, then $A$ is isomorphic to a one variable
polynomial ring over a UFD subring.
\end{lemma}

\begin{proof}
Suppose that $D\in LND(A).$ Denote by $F$ the set of fixed points of the
$G_{a}$ action on \textbf{Spec} $A$ generated by $D.$ By assumption, either
$F$ is empty, in which case $D$ has a slice \cite{F-M}, or the dimension of
$F$ is equal to one \cite{Rent}. \ In the latter case, $F$ the support of a
principal divisor $\mathfrak{D=(}f\mathfrak{)}$ for some $f\in A^{D}%
,\mathfrak{\ }$and $D(A)\subset fA$. Thus $D^{\prime}:=f^{-1}D$ is again
locally nilpotent generating a fixed point free $G_{a}$ action with a slice.
\end{proof}

Since a singular point of a factorial surface is isolated, such a surface
cannot be isomorphic to the product of a curve with a line. \ Thus

\begin{corollary}
A singular factorial surface admits no nontrivial $LND$.
\end{corollary}

The following proposition may be well known. It can be deduced from several
results in \cite{M-M} which are summarized in the proof.

\begin{proposition}
$Aut_{\mathbb{C}}R$ $\cong\ \mathbb{C}^{\ast}$ where, {for }$\lambda
\in\ \mathbb{C}^{\ast},$ $\lambda${ $(x,y,z)$ }$=(\lambda^{bc}x,\lambda
^{ac}y,\lambda^{ab}z).$
\end{proposition}

\begin{proof}
Let $\widetilde{X}$ be the quasihomogeneous factorial affine surface with
coordinate ring $R$ (whose unique singular point is the origin $0$) and
$X\equiv\widetilde{X}-\{0\}.$ Note that $Aut(X)\cong Aut(\widetilde{X}$ $)$.
\ That the mapping
\begin{align*}
G_{m}\times X  &  \rightarrow X\\
(\lambda{\ ,(x,y,z))}  &  \mapsto(\lambda^{bc}x,\lambda^{ac}y,\lambda^{ab}z)
\end{align*}
gives an action is clear. \ The quotient mapping $\pi:X\rightarrow B$,
($B\equiv X/G$) is an $A_{\ast}^{1}$ fibration, i.e. all $\pi$ fibers are
geometrically $\mathbb{C}^{\ast}$, and there are precisely three singular
fibers $F_{a},F_{b},F_{c},$of multiplicity $a,b,c$ respectively. \ In fact
$B\cong\mathbb{P}^{1},$ and any automorphism $\varphi:X\rightarrow X$
preserves the fibration, i.e. yields a group homomorphism
\[
f:Aut(X)\rightarrow Aut(\mathbb{P}^{1}).
\]
\ However, relative primeness of $a,b,c$ forces $\varphi$ to stabilize the
singular fibers and moreover $F_{a}=\pi^{-1}(\pi(F_{a})),F_{b}=\pi^{-1}%
(\pi(F_{b})),F_{c}=\pi^{-1}(\pi(F_{c})).$ Thus $\pi(F_{a}),\pi(F_{b}%
),\pi(F_{c})$ are fixed by $f(\varphi),$ and we see that $f$ is the trivial
homomorphism \cite[Cor. 4.6]{M-M}. \ Theorem 6.2 of \cite{M-M} gives the exact
sequence
\[
0\rightarrow G_{m}\rightarrow Aut(X)\rightarrow im(f)
\]
as asserted.
\end{proof}

\begin{lemma}
The restriction homomorphism {$Aut_{\mathbb{C}}A_{n,m}\rightarrow$
$Aut_{\mathbb{C}}{R}$ is surjective.}
\end{lemma}

\begin{proof}
Let $X_{n,m}$ be the affine variety with coordinate ring $A_{n,m}.$ Observe
that the mapping
\begin{align*}
G_{m}\times X_{n,m}  &  \rightarrow X_{n,m}\\
(\mu,(x,y,z,u,v))  &  \mapsto(\mu^{bc}x,\mu^{ac}y,\mu^{ab}z,\mu^{-mbc}%
u,\mu^{-nac}v)
\end{align*}
is an action inducing the $G_{m}$ action on $X$ given above.
\end{proof}

\begin{lemma}
{$\varphi\in Aut_{R}A_{n,m}$ if and only if $\varphi$ is an $R$-homomorphism
satisfying $\varphi(u,v)=(f(x,y,z)y^{n}+u,f(x,y,z)x^{m}+v)$ for some
$f\in\mathbb{C}[x,y,z]$. Consequently, $Aut_{R}A_{n,m}\cong<R,+>$ as groups. }
\end{lemma}

\begin{proof}
{ We know by corollary \ref{05.10.cor4} that $\varphi^{-1}(E)\varphi=\lambda
E$ for some $\lambda\in\mathbb{C}^{\ast}$. Define $(F,G):=\varphi(u,v)$. Also,
$\varphi(x,y,z)=(x,y,z)$. So now
\[%
\begin{array}
[c]{rl}%
(\lambda y^{n},\lambda x^{m})= & \varphi(\lambda y^{n},\lambda x^{m})\\
= & \varphi\lambda E(u,v)\\
= & \varphi(\varphi^{-1}E\varphi)(u,v)\\
= & E(F,G)\\
= & (y^{n}F_{u}+x^{m}F_{v},y^{n}G_{u}+x^{m}G_{v})
\end{array}
\]
where the subscript denotes$\ $partial derivative.}

{ Let us consider the first equation, }%
\[
{\lambda y^{n}=y^{n}F_{u}+x^{m}F_{v}.}%
\]
{ Defining $H:=F-\lambda u$, we see that $-y^{n}H_{u}=x^{m}H_{v}$. By the
following lemma \ref{05.10.diffeq} we see that $H\in R$, so }%
\[
{F=p(x,y,z)+\lambda u.}%
\]
{ The second equation yields $\lambda x^{m}=y^{n}G_{u}+x^{m}G_{v}$. Defining
$H:=G-\lambda v$, yields $-x^{m}H_{v}=y^{n}H_{u}$, which by the following
lemma \ref{05.10.diffeq} yields $H=q(x,y,z)$ and thus $G=q(x,y,z)+\lambda
v$. Now}%

\[%
\begin{tabular}
[c]{lll}%
$0$ & $=$ & ${\varphi(x^{m}u-y^{n}v-1)}$\\
& $=$ & ${x^{m}\varphi(u)-y^{n}\varphi(v)-1}$\\
& $=$ & ${x^{m}F-y^{n}G-1}$\\
& $=$ & ${x^{m}(p+\lambda u)-y^{n}(q+\lambda v)-1}$\\
& $=$ & ${x^{m}p-y^{n}q+\lambda(x^{m}u-y^{n}v)-1}$\\
& $=$ & ${x^{m}p-y^{n}q+\lambda-1.}$%
\end{tabular}
\
\]

{ Thus $\lambda=1$ and $p=y^{n}f(x,y,z)$ and $q=x^{m}f(x,y,z)$ for some $f$.
It is not difficult to check that the constructed objects are well-defined
homomorphisms which are isomorphisms. }
\end{proof}

\begin{lemma}
{ \label{05.10.diffeq} If $H\in A_{n,m}$ such that $-y^{n}H_{u}=x^{m}H_{v}$,
then $H\in R$. }
\end{lemma}

\begin{proof}
{ We can find polynomials $p_{i}(v)\in$}${R[v]=}${ $\mathbb{C}$}$[x,y,z][v]${
such that $H=\sum_{i=0}^{d}p_{i}u^{i}$ for some $d\in\mathbb{N}$. Requiring
$deg_{z}(p_{i})<c$ for each $i\in\mathbb{N}^{\ast}$, and $deg_{x}(p_{i})<m$
for each $i\in\mathbb{N}^{\ast},i\not =1$, then the $p_{i}$ are unique
(because of the equality $x^{m}u=y^{n}v+1$ and $z^{c}=-x^{a}-y^{b}$). The
equation $-y^{n}H_{u}=x^{m}H_{v}$ yields
\[
\sum_{i=0}^{d-1}-(i+1)y^{n}p_{i+1}u^{i}=\sum_{i=0}^{d}x^{m}p_{i,v}u^{i}%
\]
where }$p_{i,v}\equiv\frac{\partial p_{i}}{\partial v}.$ \ {Substitute
$y^{n}v+1$ for $x^{m}u$ to obtain a unique representation:}%
\[%
\begin{array}
[c]{rl}%
\sum_{i=0}^{d-1}-(i+1)y^{n}p_{i+1}u^{i}= & x^{m}p_{0,v}+\sum_{i=0}^{d-1}%
(y^{n}v+1)p_{i+1,v}u^{i},
\end{array}
\]
so
\[%
\begin{array}
[c]{rl}%
-y^{n}p_{1}= & x^{m}p_{0,v}+(y^{n}v+1)p_{1,v}%
\end{array}
\]
and%
\[%
\begin{array}
[c]{rl}%
-(i+1)y^{n}p_{i+1}= & (y^{n}v+1)p_{i+1,v}%
\end{array}
\]
for each $i\geq1$.

{ Let $i\geq1$ and assume that $p_{i+1}$ has degree $k$ with respect to $v$.
Let $\alpha(x,y,z)$ be the top coefficient of $p_{i+1}$, seen as a polynomial
in $v$. Then $-(i+1)y^{n}\alpha=y^{n}k\alpha$, but that gives a contradiction.
So for each $i\geq1:p_{i+1}=0$. This leaves the equation $0=x^{m}p_{0,v}$
which means that $p_{0}\in\mathbb{C[}x,y,z]$. Thus $H=p_{0}u^{0}\in
\mathbb{C[}x,y,z]$. }
\end{proof}

{ We conclude this section with a statement of the theorem just proved:}

\begin{theorem}
{$Aut_{\mathbb{C}}A_{n,m}$ is generated by the maps\newline}
\end{theorem}

\begin{enumerate}
\item $(x,y,z,u,v)\mapsto{ (x,y,z,f(x,y,z)y^{n}+u,f(x,y,z)x^{m}+v)}$ { for
}$f\in R,$

\item $(x,y,z,u,v)\mapsto(\mu^{bc}x,\mu^{ac}y,\mu^{ab}z,\mu^{-mbc}u,\mu
^{-nac}v)$ for $\lambda\in\mathbb{C}^{\ast}${. }
\end{enumerate}

{Thus $Aut_{\mathbb{C}}A_{n,m}\cong\mathbb{C}^{\ast}$}${\ltimes}<R,+>${ . }

Note that {$Aut_{\mathbb{C}}A_{n,m}$ is nonabelian.}

\section{Examples}

\begin{example}
Let $R={\mathbb{C[}X,Y,Z]/(X^{a}+Y^{b}+Z^{c})}$ where $a,b,c\ $are pairwise
relatively prime positive integers satisfying $\frac{1}{a}+\frac{1}{b}%
+\frac{1}{c}<1{.}${ \ Then }$A_{n,m}\times\mathbb{C}\cong A_{n^{\prime
},m^{\prime}}\times\mathbb{C}$ for all $(n,m),(n^{\prime},m^{\prime})$ but
$A_{n,m}\cong A_{n^{\prime},m^{\prime}}$ if and only if $(n,m)=(n^{\prime
},m^{\prime}).$ \ Hence the $X_{n,m}\equiv\mathbf{Spec}A_{n,m}$ are the
desired counterexamples to the generalized affine cancellation problem.
\end{example}

\begin{proof}
Since the $\mathbf{Spec}A_{n,m}$ are all total spaces for principal $G_{a}$
bundles over $\mathbf{Spec}R-\{(0,0)\},$ the first assertion is clear. \ Write
$A_{n,m}=R[u,v]$ where $x^{m}u-y^{n}v=1,$ and $A_{n^{\prime},m^{\prime}%
}=R[u^{\prime},v^{\prime}]$ where $x^{m^{\prime}}u^{\prime}-y^{n^{\prime}%
}v^{\prime}=1.$Since $\frac{1}{a}+\frac{1}{b}+\frac{1}{c}<1$, $ML(A_{n,m})=R$
and an isomorphism $\Phi:A_{n,m}\cong A_{n^{\prime},m^{\prime}}$ will restrict
to an automorphism of $R.$ Thus, possibly after a composition with an
automorphism of $R$,
\[
\Phi(x)=x,\text{ }\Phi(y)=y,\text{ }\Phi(z)=z.
\]
\ Let $D$ $\in$ $LND(A_{n,m})$ (resp. $D^{\prime}\in LND(A_{n^{\prime
},m^{\prime}})$) satisfy
\begin{align*}
D  &  :v\mapsto x^{m}\mapsto0,u\mapsto y^{n}\mapsto0\\
D^{\prime}  &  :v^{\prime}\mapsto x^{m^{\prime}}\mapsto0,u^{\prime}\mapsto
y^{n^{\prime}}\mapsto0.
\end{align*}
Since $LND(A_{n,m})=RD$ and $D,D^{\prime}$ are irreducible derivations, the
locally nilpotent derivation $\Phi^{-1}D^{\prime}\Phi=rD$ for some $r\in
R^{\ast}=\mathbb{C}^{\ast}.$

Set $K=qf(R)$, identify $K\otimes_{R}A_{n,m}$ with $K[v]$, $K\otimes
_{R}A_{n^{\prime},m^{\prime}}=K[v^{\prime}]$, and note that $K[\Phi
(v)]=K[v^{\prime}].$ \ Thus
\[
\Phi(v)=\alpha v^{\prime}+\beta\text{ for some }\alpha,\beta\in K.
\]
A calculation reveals that
\[
\Phi^{-1}D^{\prime}\Phi(v)=\Phi^{-1}(\alpha)x^{m^{\prime}}=rx^{m}.
\]
so that $\alpha x^{m^{\prime}}=rx^{m}.$

We obtain%
\[
\Phi(v)=x^{m-m^{\prime}}v^{\prime}+\beta
\]
from which we conclude that $D^{\prime2}(\Phi(v))=0.$ \ A symmetric argument
yields that $D^{\prime2}(\Phi(u))=0$ as well. \ Thus
\begin{align*}
\Phi(u)  &  =r_{1}u^{\prime}+r_{2}v^{\prime}+r_{3}\\
\Phi(v)  &  =s_{1}u^{\prime}+s_{2}v^{\prime}+s_{3}%
\end{align*}
with $r_{i},s_{j}\in R,$ and $r_{1}s_{2}-r_{2}s_{1}\in R^{\ast}.$

If $m>m^{\prime}$, then $\beta\in K\cap A_{n^{\prime},m^{\prime}}=R$, so that
$s_{1}\in xR,$ $s_{2}=\mu^{\prime}x^{m-m^{\prime}},$ and $s_{3}=\beta.$ But in
this case
\[
r_{1}s_{2}-r_{2}s_{1}\in xR\nsubseteqq R^{\ast}.
\]
Thus $m\leq m^{\prime}$, but the identical argument with the roles of $\Phi$
and $\Phi^{-1}$reversed will show $m=m^{\prime},$ and the symmetric argument
with the roles of $u$ and $v$ reversed will show $n=n^{\prime}.$
\end{proof}

\end{document}